\documentclass[10pt]{article}

\usepackage{amsmath}
\usepackage{amssymb}
\usepackage{bm}
\usepackage{color}
\usepackage{graphicx}
\usepackage{hyperref}
\usepackage{multirow}
\usepackage[square, numbers, comma, sort&compress]{natbib}  
\usepackage{vector}  
\usepackage{verbatim} 

\newtheorem{prop}{\indent{\sc Proposition}}
\newtheorem{rmk}{\indent{\sc Remark}} 
\newtheorem{thm}{\indent{\sc Theorem}}

\newcommand{\benu}{\begin{enumerate}}
\newcommand{\eenu}{\end{enumerate}}

\newcommand{\beit}{\begin{itemize}}
\newcommand{\eeit}{\end{itemize}}

\newcommand{\be}{\begin{eqnarray}}
\newcommand{\bec}{\begin{center}}
\newcommand{\eec}{\end{center}}
\newcommand{\beo}{\begin{eqnarray*}}
\newcommand{\Bl}{\Bigl(}
\newcommand{\Br}{\Bigr)}

\newcommand{\ee}{\end{eqnarray}}
\newcommand{\eeo}{\end{eqnarray*}}

\newcommand{\fin}{\frac{i}{n}}
\newcommand{\fini}{\frac{i+1}{n}}

\newcommand{\fracn}{\frac{1}{\sqrt{n}}}

\newcommand{\fpi}{\frac{1}{\sqrt{2 \pi}}}

\newcommand{\intt}{\int \limits_{0}^{t}}

\newcommand{\mA}{\mathcal{A}}
\newcommand{\mb}{\mathbf{b}}
\newcommand{\mB}{\mathcal{B}}

\newcommand{\mepsilon}{\mathbf{\epsilon}}

\newcommand{\mbR}{\mathbb{R}}

\newcommand{\mR}{\mathbb{R}}

\newcommand{\mx}{\mathbf{x}}
\newcommand{\mX}{\mathbf{X}}
\newcommand{\mY}{\mathbf{Y}}
\newcommand{\my}{\mathbf{y}}
\newcommand{\mz}{\mathbf{z}}
\newcommand{\mZ}{\mathbf{Z}}

\newcommand{\n}{\nonumber}

\newcommand{\R}{\Rightarrow}

\newcommand{\tg}{\textcolor{black}}
\newcommand{\tr}{\textcolor{black}}

\newcommand{\mcX}{\mathcal{X}}

\begin{document}

\title{Langevin type limiting processes for Adaptive MCMC}
\author{ G K Basak \footnote{Stat-Math Unit, Indian Statistical Institute, 203 B T Road, Kolkata 108, India.} \& Arunangshu Biswas 
\footnote{Communicating author, Dept. of Statistics, Presidency University, Kolkata 73, India. Email: arunanangshu.stats@presiuniv.ac.in}}
\date{ }
\maketitle

\begin{abstract}
Adaptive Markov Chain Monte Carlo (AMCMC) is a class of MCMC algorithms where the proposal distribution changes at every iteration of the chain. In this case it is important to verify that such a Markov Chain indeed has a stationary distribution. In this paper we discuss a diffusion approximation to a discrete time AMCMC. This diffusion approximation is different when compared to the diffusion approximation as in Gelman, Gilks and Roberts (1997) \cite{Gelman} where the state space increases in dimension to $\infty$. In our approach the time parameter is sped up in such a way that the limiting distribution (as the mesh size goes to 0) is to a non-trivial continuous time diffusion process. 

\end{abstract}


\section{Introduction}
\label{introduction}
Markov Chain Monte Carlo (MCMC) methods allow us to generate samples from an arbitrary distribution $\psi(\cdot)$ known up to a scaling factor. The algorithm, consists in sampling a Markov chain $\{X_k, k \ge 0 \}$ on a state space $\mcX$ with transition probability $P$ admitting $\psi(\cdot)$ as its unique invariant distribution.\\
In most MCMC algorithms known so far, the transition probability $P$ depends on some tuning parameter $\theta$ defined on some space $\Theta$ which is possibly infinite dimensional. The success or failure of the algorithm depends on the choice of $\theta.$ To see this, consider the Metropolis-Hastings (MH) algorithm  where we assume that the target distribution $\psi(\cdot)$ admits a density w.r.t to the Lebesgue measure, which we also denote as $\psi(\cdot)$ on $\mcX = \mbR.$. Then given that the chain is at $x \in \mbR$ we choose a candidate $y$ according to the proposal transition density $q(x,\cdot)$ and is accepted with the probability 
$\alpha(x,y) = \min \{1, \frac{\psi(y)}{\psi(x)} \frac{q(y,x)}{q(x,y)}\}$. 
A common choice for $q(x, \cdot)$ is the Normal density function with mean $x$ and variance $\theta^2 < \infty$. Since in that case $q(x,y) = q(y,x), \ \alpha$   takes the simple form 
\be 
\alpha(x,y) &=& \min \{1, \frac{\psi(y)}{\psi(x)} \}. \label{MH accept}
\ee Therefore the algorithm would proceed as follows: 1) Start with an initial $X_0$, \ 2) For any $n \ge 1$ generate a sample $y$ from $q(X_n, \cdot)$ and accept that with probability $\alpha(X_n,y).$ Call it $X_{n+1}$. It can be shown that this $\{X_n \}$ is a Markov Chain on the support of $q(x, \cdot)$ and is invariant with respect to $\pi(\cdot).$\\
The main drawback of the algorithm is that the speed of convergence of the Markov chain to the invariant distribution depends on the choice of $q(\cdot, \cdot)$. Bad choice of the proposal distribution makes the convergence to stationary too slow. \\
The problem of the optimal choice of the proposal distribution $q(\cdot, \cdot)$ was dealt in a paper by Gelman, Gilks and Roberts  \cite{Gelman}. In their paper the target was a  $d$-dimensional with i.i.d components. The proposal is multivariate normal with mean zero and dispersion $\frac{\sigma^2}{d} \textbf{I}_d$. Their interest was the infinite dimensional process as $d \to \infty$. Suitably scaling the time and space parameter the limiting continuous time process was obtained. It turned out that the acceptance rate optimizing the efficiency of the process as $d \to \infty$ converges to $0.234$.\\
Adaptive MCMC, introduced by Haario, Saksman and Tamminen (2001) \citep{Haario}, are a class of algorithms that adjusts the transition kernel according to the previous values of the chain. In this case since the transition kernel change at each iteration it convergence to stationarity is not automatically guaranteed. Sufficient conditions are given in Roberts and Rosenthal (2005), \cite{Roberts}. 

In this paper we obtain the invariant distribution of a suitably defined AMCMC, after performing the diffusion approximation procedure to the process. For details of the diffusion approximation see, for example \citep{Nelson}. Our choice of the AMCMC arises from the fact that the adaptation parameter (also called \textit{tuning/ scaling parameter}) should depend on whether the sample generated from the proposal distribution is accepted or not. If accepted, then the scaling parameter should increase by some amount and if not, the scaling parameter should decrease. 
The outline of the paper is as follows: In Section \ref{defn} we give the definition of the proposed AMCMC (a partial variation of this algorithm was suggested by Prof. P. Green in a personal communication.) In Section \ref{diffusion approximation} we give the details of the diffusion approximation procedure. Section \ref{Theorems} gives the main Theorem of the paper. Section \ref{drift and diffusion coefficients} deals with some computations required in the proof of the Theorem in Section \ref{Theorems}. In Section \ref{multsec} we discuss the case where the target distribution is a general multivariate distribution and the proposal is multivariate Normal ($\mathbf{0}, \mathbf{I}_p$). We describe this procedure in details since the proof depends on the elliptical symmetric property of the multivariate Normal distribution  whereas the univariate case relies only on the symmetric and finite second momemts. Some simulations are provided in Section \ref{sim}. We end with some concluding remarks in Section \ref{summary}.

\section {Definition of the Adaptive MCMC algorithm}
\label{defn}

{We assume that the target distribution $\psi(\cdot)$ is univariate and $\frac{\psi'(x)}{\psi(x)}$ grows linearly in $x$. (The reason for this choice is explained in Remark \ref{gronwall}). 

\textbf{Algorithm 1:}
\begin{enumerate}

\item Select arbitrary $\{X_0,\theta_0,\xi_0 \}\in \mathbb{R}\times( 0,\infty)\times\{0,1\}$
 where $\mathbb{R}$ is the state space which may be the real line or an interval of the same. Set $n=1$.

\item
Propose a new move, say Y, where \ \ \
$Y \sim N(X_{n-1}, \theta_{n-1}$).

\item
Accept the new point with probability $\alpha(X_{n-1}, Y)= \min\{1, \frac{\psi(Y)}{\psi(X_{n-1})}\}$.\\
If the point is accepted, set $X_{n}=Y ,\ \xi_n=1$; else $X_{n}=X_{n-1},\ \xi_n=0$.

\item
$\theta_{n}=\theta_{n-1} e^{\frac{1}{\sqrt{n}}(\xi_n - q)}, \ \ q>0, 
\ \ \Leftrightarrow \log(\theta_n) = \log(\theta_{n-1}) + \frac{1}{\sqrt{n}}(\xi_n - q), \ \ q>0$.

\item
$n \leftarrow n+1$,
and go to step 2.

\end{enumerate}

{The above algorithm is equivalent to the following:} \\ \ \ \
\textbf{Algorithm 1$'$:} \ \ \
\begin{enumerate}
\item Select arbitrary $\{X_0,\theta_0,\xi_0 \}\in \mathbb{R}\times(0,\infty)\times\{0,1\}, $
 where $\mathbb{R}$ is the state space. Set $n=1$.
\item Given  $X_{n-1}, \theta_{n-1}, \epsilon_{n-1} $ generate $$ \xi_n \sim Bernoulli \Bl \min \Bl 1, \frac{\psi(X_{n-1} + \theta_{n-1} \epsilon_{n-1}) }{\psi(X_{n-1})} \Br \Br$$
and then
\be 
X_{n} &=& X_{n-1} + \theta_{n-1} \xi_{n} \epsilon_{n-1} \label{stoc. app.}
\ee
where $\epsilon_{n-1} \sim N(0,1)$, 
\item $\theta_{n}=\theta_{n-1} e^{\frac{1}{\sqrt{n}}(\xi_n - q)}, \ \ q>0,  \ \ \Leftrightarrow \log(\theta_n) =\log(\theta_{n-1})+ \frac{1}{\sqrt{n}}(\xi_n - q), \ \  q>0.$
\item
$n \leftarrow n+1$
and go to step 2.
 \end{enumerate}

Let us describe the algorithm. $\theta_n$ is the proposal scaling (tuning) parameter which is adaptively tuned 
depending on 
whether the previous sample was accepted \tr{or rejected}. If the sample was accepted then the proposal variance will increase 
allowing the chain to explore more regions in the state space. If the past sample was rejected then the 
variance will decrease making the move a more conservative one. Here $q$  is a benchmark; for multivariate Normal target density, where the components are independent, the value 0.238 is often appropriate, see Gelman \textit{et al}. \citep{Gelman}. For \tr{a further} generalization see Bedard \citep{Bedard2}. {The tuning parameter can also be made to be dependent not only on whether the previous sample was accepted but also on the proportions of \tr{samples} accepted in the history of the chain. However, that is not done in this paper.} \\
Now, 
\beo  
E(\theta_n \xi_{n+1} \epsilon_n| \theta_n, X_n) &=& E(E(\theta_n \xi_{n+1} \epsilon_n| \theta_n, X_n, \epsilon_n)| \theta_n, X_n) \\
&=& E (\theta_n \epsilon_n E (\xi_{n+1}|X_n, \theta_n, \epsilon_n)| \theta_n, X_n) 
= E (\theta_n \epsilon_n \mathbb{P} (\xi_{n+1} = 1 |X_n, \theta_n, \epsilon_n)| \theta_n, X_n) \\
&=& E(\theta_n \epsilon_n \min \{1, \frac{\phi(X_n + \theta_n \epsilon_n)}{\phi(X_n)}\}| \theta_n, X_n)\\
&=& \theta_n \int_{\mathbb{R}}  \epsilon \min \{1, \frac{\phi(X_n + \theta_n \epsilon)}{\phi(X_n)}\} \phi(\epsilon) d \epsilon  \\
&:=& p(X_n, \theta_n).
\eeo
Therefore Equation (\ref{stoc. app.}) can be written as 
\beo
\Delta X_n = X_{n+1} -  X_n &=& p(X_n, \theta_n) + \Bl \theta_n \xi_{n+1} \epsilon_n - p(X_n, \theta_n)\Br.
\eeo
Define $M_n = \theta_n \xi_{n+1} \epsilon_n - p(X_n, \theta_n)$, then $E(M_n| X_n, \theta_n) = 0,$ which means that $\{M_n\}$ is a martingale difference sequence (w.r.t its natural filtration). This bears similarity with the Stochastic Approximation procedure which was introduced by Robbins and Monro \citep{Monro}.
For a recent review see Borkar \citep{Borkar} and references therein.

Here we embed the discrete time chain into a continuous time stochastic process. 
This technique has been applied \tr{to} diverse fields, for example, econometric modelling (Nelson  \citep{Nelson}), branching processes (Ethier and Kurtz \citep{Ethier}). \tr{One advantage is that we can apply standard tools in continuous time stochastic processes, which are not available for discrete time AMCMC.}
The next section gives details of the diffusion approximation technique.

\section{Diffusion Approximation}
\label{diffusion approximation}
In this section we first  present conditions developed by Stroock and Varadhan \citep{Stroock} for a sequence of stochastic processes satisfying a stochastic difference equations to converge weakly to an It\^o Process. \\
\tr{ Here is the formal set up: Let $D([0, \infty), \mR^n)$ be the space of mappings from $[0, \infty)$ into $\mR^n$ that are continuous from the right with left limits and let $\mB(\mR^n)$  denote Borel sets in $\mR^n$. $D$ is a metric space when endowed with the Skorokhod metric (see Billingsley \citep{Billingsley}). For each $h>0$ let $\mathcal{M}_{kh}$ be the $\sigma$-algebra generated by the random variables $\mX_{0,h}, \mX_{h,h}, \ldots, \mX_{kh,h}$ for $k \ge 1$ and let $\nu_h$ be a probability measure on $(\mR^n, \mB(\mR^n))$. For each $h>0$ and each $k \ge 1$ let $\Pi_{kh,h}$ be a transition kernel for a homogeneous Markov chain i.e.,}
\benu
\item $\Pi_{kh,h}(\mx, \cdot)$ is a probability measure on $(\mR^n, \mB(\mR^n))$ for all $\mx \in \mR^n;$
\item $\Pi_{kh,h}(\cdot, A)$ is a $\mB(\mR^n)$ measurable for all $A \in \mB(\mR^n).$
\eenu
For each $h>0$, let $P_h$ be the probability measure on $D([0, \infty), \mR^n)$ such that
\be
P_h\Bl \mX_{0,h}  \in A \Br  &=& \nu_h(A) \ \ \forall  \ A \in \mB(\mR^n), \label{condn 1}\\
P_h \Bl \mX_{t,h} = \mX_{kh,h}, \ kh \le t \le (k+1)h \Br &=& 1, \label{condn 2} \\
P_h \Bl \mX_{(k+1)h, h} \in A| \mathcal{M}_{kh}\Br &=& \Pi_{kh,h}(\mX_{kh,h}, A).  \label{condn 3}
\ee
{almost surely under $P_h \ \forall \ k \ge 0$ and $A \in \mB(\mR^n)$ }.\\
For each $h>0$, equation (\ref{condn 1}) specifies the distribution of the random starting point.  In Equation  (\ref{condn 2}) we construct a continuous time process from the discrete time process by making $X_{t,h}$ a step function with jumps at time $h, 2h, 3h, \ldots$ etc.  \tr{Equation (\ref{condn 3}) states that for a fixed $h > 0 $, $\{X_{kh,h}, k \ge 1 \}$ is a Markov Chain with $\Pi_{kh,h}(\cdot, \cdot)$ as the transition kernel.} \\
We next define the \tg{infinitesimal} diffusion and drift coefficients \tr{for any $t, h >0$} as :
\be
a_h(\mx,t) &=& h^{-1} \int_{\mR^n} (\my-\mx)(\my-\mx)' \Pi_{[t/h]h,h}(\mx, d\my) \nonumber \\
&=& h^{-1} D(\mX_{(k+1)h,h} | \mX_{kh,h} = \mx) \ \mbox{for any $k \ge 1$}; \n \\
b_h(\mx,t) &=& h^{-1} \int_{\mR^n} (\my-\mx) \Pi_{[t/h]h,h}(\mx,d\my) = h^{-1} E(\mX_{(k+1)h,h} - \mx |\mX_{kh,h} = \mx) \ \mbox{for any $k \ge 1$};\n \\
\Delta_{h, \epsilon}(\mx,t) &=& h^{-1} \int_{|| \my - \mx || > \epsilon} \Pi_{[t/h]h,h}(\mx,d\my) \n \\
&=&  h^{-1} P(|| \mX_{(k+1)h,h} - \mX_{kh,h} || > \epsilon \ | \ \mX_{kh,h} = \mx)\ \mbox{for any $k \ge 1$}, \label{drift equations}
\ee
\tr{where $D(\mX_{(k+1)h,h} |\mX_{kh,h} = \mx) $ and $\ E(\mX_{(k+1)h,h} - \mx |\mX_{kh,h} = \mx)$ are the conditional dispersion and conditional expected deviation given that the value of $\mX_{kh,h}$ is $\mx$ respectively. $a_h(\mx,t)$ and $b_h(\mx,t)$ are measures of the second moment and drift per unit of time respectively. $\Delta_{h, \epsilon}(\mx, t)$ is the conditional probability  of a jump of size $\epsilon$ or greater per unit of time. 
The convergence results that we present below will require that $a_h(\mx,t)$ and $b_h(\mx,t)$ converge to a finite limits and $\Delta_{h, \epsilon}(\mx,t)$ goes to zero for all $\epsilon > 0$ as $h \downarrow 0$}. In particular we assume the following, see \citep{Stroock}:\\ 
\textbf{Assumptions}
\benu
\item There exists a locally bounded measurable mapping $a(\mx,t): \mR^n \times [0, \infty) \to M_{n \times n}^+$ which are continuous in $x$ for each $t \ge 0$, and \tr{$b(\mx,t):\mathbb{R}^n \times [0, \infty) \to \mR^n$} such that: 
\beo
\lim_{h \downarrow 0} ||a_h(\mx,t) - a(\mx,t)||  &=& 0; \\
\lim_{h \downarrow 0} ||b_h(\mx,t) - b(\mx,t)||  &=& 0; \\
\lim_{h \downarrow 0} \Delta_{h, \epsilon}(\mx,t) &=&0,
\eeo
where $M_{n \times n}^+$ denotes the space of all $n \times n$ non-negative definite matrices and $|| \cdot||$ is the matrix/vector  norm defined as:

\tr{
$$
|| A || =\left\{
\begin{array}{ll}
{[A^T A]}^\frac{1}{2} & \mbox{if $A$ is a column vector}\\
{[trace (A^T A)]}^\frac{1}{2} & \mbox{if $A$ is a matrix}.
\end{array}
\right.
$$
}
\item There exists a locally bounded measurable mapping $\sigma(\mx,t)$ form $\mR^n \times [0, \infty) \to M_{n \times n}$ which are continuous in $x$ for each $t\ge 0$, such that for all $\mx \in \mR^n$ and all $t \ge 0$, 
$$ a(\mx, t) = \sigma(\mx,t) \sigma(\mx,t)',$$
where $M_{ n \times n}$ denotes the space of all $n \times n$ matrix.
\item As $h \downarrow 0, \ X_{0,h}$ converges in distribution to a random variable $X_0$ with a probability measure $\nu_0$  on $(\mR^n, \mB(\mR^n))$;
\item $\nu_0, a(\mx,t)$ and $b(\mx,t)$ uniquely specify the distribution of a diffusion process $\mX_t$, with the initial distribution $\nu_0$, the the diffusion matrix $a(\mx,t)$ and the drift vector $b(\mx,t)$.
\eenu
Under the assumption we have the following Proposition. For a proof see Stroock and Varadhan \citep{Stroock}.
\begin{prop}
Under Assumptions 1 - 4, the sequence of $\mX_{h,t}$ process defined by Equations (\ref{condn 1}) - (\ref{condn 3}) converges weakly (i.e., in distribution) as $h \downarrow 0$ to the $\mX_t$ process defined by the stochastic integral equation 
{
\be
\mathbf{X}_t = \mathbf{X}_0 + \int_{0}^{t} b(\mathbf{X}_s,s)ds + \int_{0}^{t} \sigma(\mathbf{X}_s, s) d W_{n,s} \label{basic}
\ee
}
\end{prop}
where $W_{n,t}$ is an $n$-dimensional standard Brownian motion, independent of $\mX_0$ and where for any $A \in \mB(\mR^n), P(X_0 \in A) = \nu_0(A)$. Such an $\mX_t$ process exists and is unique upto a distribution. 

Next we  embed the discrete time Algorithm 1$'$ defined in Section \ref{defn} in a continuous time process that has decreasing step sizes. For fixed $n \ge 1,$ we partition the half line $[0,\infty)$ into sub intervals of length $\frac{1}{n}$. We start with the fixed point $x_0$. Now given the value of the process at time $\frac{i}{n}$, i.e., $X_n \Bl \frac{i}{n}\Br = x$, we propose a value following the $N \Bl x, \frac{1}{\sqrt{n}}\theta_n\Bl \frac{i}{n} \Br \Br$ distribution. We have the correction factor $\frac{1}{n}$ multiplied with the variance to incorporate the diminishing adaptation condition, so that the difference between the proposal kernel at times $\frac{i}{n}$ and $\frac{i+1}{n}$ goes to zero as $n \to \infty$. This proposed value is accepted with the usual MH acceptance probability given in  (\ref{MH accept}) at time $\frac{i+1}{n}$. The indicator variable denoting whether the proposed value is  accepted is denoted by $\xi_n\Bl \frac{i}{n}\Br$. Similar approximation is done with the tuning parameter $\theta_n(\cdot)$ starting with the initial value $\theta_0$.

\subsection{{Embedding in continuous time of discrete AMCMC}}
\label{Embedding in continuous time of discrete AMCMC}
The following gives the embedding of the discrete AMCMC into continuous times
{state variable $X_n(\cdot)$}
\be
X_n(0)&=&x_0 \in \mathbb{R}; \n \\
X_n\Bl \frac{i+1}{n} \Br &=& X_n\Bl \frac{i}{n} \Br + \frac{1}{\sqrt{n}} \theta_{n}\Bl\frac{i}{n}\Br \xi_n\Bl\frac{i+1}{n}\Br \epsilon_n\Bl\frac{i+1}{n}\Br, \ \ \mbox{i=0, 1, \ldots},\n \\
X_n(t)&=& X_n\Bl{\frac{i}{n}}\Br ,  \ \ \mbox{if $\frac{i}{n} \le t < \frac{i+1}{n}$} \ \ \mbox{for some integer $i$.} \label{X}
\ee
Here, $\xi_n(\frac{i+1}{n})$ conditionally follows the Bernoulli distribution given by:
\beo
& & P\Bl \xi_n(\frac{i+1}{n})=1|X_n(\frac{i}{n}),\  \theta_n\Bl\frac{i}{n}\Br, \ \epsilon_n \Bl \fini \ \Br \Br \n \\
&=& \min \Bigl \{\frac{\psi(X_n\Bl\fin\Br+ \frac{1}{\sqrt{n}}\theta_n\Bl\frac{i}{n}\Br \epsilon_n\Bl \fini \Br\Br}{\psi \Bl X_n \Bl \frac{i}{n}\Br \Br},1 \Bigr\}
\eeo
and $\{\epsilon_n(\frac{i}{n}), i \ge 1 \ \} $ are all independent $N(0,1)$ random variables. This \tr{distribution} of $\xi_n(\cdot)$ comes directly from the form of the MH acceptance probability given in (\ref{MH accept}).

\textbf{Tuning parameter \textbf{$\theta_n(\cdot)$}}\\
The $n^{th}$ approximation to the tuning parameter $\theta(\cdot)$ is defined as :
\be
\theta_n(0) &=& \theta_0 \in \mathbb{R^+}, \n \\
\theta_n\left({\frac{i+1}{n}}\right) &=& \theta_n\left({\frac{i}{n}}\right) e^{\frac{1}{\sqrt{n}} (\xi_n(\frac{i+1}{n})-q_n({\frac{i}{n}}))}, \ \ \mbox{i=0, 1, \ldots},\n \\
\mbox{and \ } \theta_n(t)&=& \theta_n(\frac{i}{n}) , \ \ \mbox{if $\frac{i}{n} \le t < \frac{i+1}{n}$ for some integer $i$}. \label{theta}
\ee

\tr{In the original discrete AMCMC the benchmark value of $q$, given in Step 3 of Algorithm 1$'$, was kept fixed. However if that is also done in the continuous AMCMC in Equation  (\ref{theta}) then the tuning parameter $\theta_n$ will converge to 
\bec 
$
\left \{
\begin{array}{ll}
\infty & \mbox{if} \ q < 1; \\
0      & \mbox{if} \ q = 1.
\end{array}
\right. 
$
\eec
It is exactly for this reason the constants in the tuning parameter given in Equation  (\ref{theta}) is an increasing function of n (also depending on a constant $q>0$) that converges to 1 as $n \to \infty$. In particular, for our example, we have $q_n\Bl\frac{i}{n}\Br = 1 - \frac{q}{\sqrt{n}}$ for some $q >0$.}

For comparison purposes we also embed the discrete time standard MCMC (SMCMC) in continuous times. The SMCMC algorithm is almost similar to the AMCMC, except for the fact that the tuning parameter given by $\theta(\frac{i}{n})$ corresponding to SMCMC is kept fixed at a constant level $\theta_0$, that is unchanged in the iterations. This is done in the next subsection .

\subsection{{Embedding in continuous times of SMCMC}}
The continuous time process corresponding to SMCMC will therefore be :\\
\be
X_n(0)&=&x_0 \in \mathbb{R}; \n \\
X_n\Bl \frac{i+1}{n} \Br &=& X_n\Bl \frac{i}{n} \Br + \frac{1}{\sqrt{n}} \theta_0 \xi_n\Bl\frac{i+1}{n}\Br \epsilon_n\Bl\frac{i+1}{n}\Br, \ \ \mbox{i=0, 1, \ldots}, \ \theta_0 \in \mR^+=(0,\infty), \n \\
X_n(t)&=& X_n\Bl{\frac{i}{n}}\Br ,  \ \ \mbox{if $\frac{i}{n} \le t < \frac{i+1}{n}$} \ \ \mbox{for some integer $i$.} \label{XSMC}
\ee
where $\xi_n \Bl  \frac{i}{n}\Br$ has the same conditional distribution with $\theta_n$ replaced by $\theta_0$ where $\theta_0$ is the fixed constant that is not updated in the iterations.\\
The following main Theorem of this paper tells the outcome of the diffusion approximation of the Discrete AMCMC defined through Equations (\ref{X}) to (\ref{theta}) and that of the SMCMC defined through (\ref{XSMC}).

\section{Main Theorem}
\label{Theorems}
\setcounter{equation}{0}
\begin{thm}
\label{main thm 2}
\benu
\item
 $\mathbf{Y}_n(t) := \Big(X_n(t), \ \  \theta_n(t) \Big)$ (where $X_n(t)$ and $\theta_n(t)$ is given by (\ref{X}) and (\ref{theta}) respectively)  converges weakly to a diffusion process \tg{which} is the solution to the SDE, 
\be 
d \mathbf{Y}_t &=& b(\mathbf{Y}_t)dt+ \sigma(\mathbf{Y}_t)d\mathbf{W}_t. \label{coupled system}
\ee
Here,
\beo
b(\mathbf{Y}_t)=\left(  \frac{\theta_t^2}{2} \frac{\psi'(X_t)}{\psi(X_t)} , \ \  \theta_t \left(q - \frac{ \theta_t}{\sqrt{2 \pi}} \frac{|\psi'(X_t)|}{\psi(X_t)} \right) \right)^{T},\\
\eeo
\\
and
\beo
\sigma(\mathbf{Y_t})=
\left(
\begin{array}{cc}
\theta_t & 0 \\
0 & 0
\end{array}
\right),
\eeo
\item Similarly the SMCMC converges weakly to a diffusion to the process which is the solution to the SDE
\be
  d X_t = \frac{\psi'(X_t)}{\psi(X_t)} \frac{\theta_0^2}{2} dt + \theta_0 dW_t \label{X MC}.
  \ee

\eenu
and $\mathbf{W}_t$ is a two dimensional Brownian motion. See Remarks \ref{conditions} for more details on the conditions on $\psi(\cdot)$.
Here $x^T$ is the transpose of a vector (or, a matrix) $x$.
\end{thm}

\textbf{Proof.}
Firstly, note that since $\mathbf{Y}_n(\frac{i}{n}) := \left(X_n(\frac{i}{n}), \ \ \theta_n(\frac{i}{n})  \right)$ 
is a homogeneous Markov chain it defines a transition kernel 
\beo
\Pi_n(\mathbf{y},A) = P\left( Y_n(\frac{i+1}{n}) \in A|\mathbf{Y}_n(\frac{i}{n})=\mathbf{y} \right), \ \ \forall \my \in \mathbb{R} \times \mathbb{R^+} \ \ \mbox {and} \ \ \forall A \in \mathcal{B}(\mathbb{R} \times \mathbb{R}^{+}) .
\eeo
Note that since the initial points of the AMCMC and the SMCMC is fixed at $(x_0, \theta_0)$ Assumption 3 of Section \ref{diffusion approximation} is automatically satisfied, where $\nu_0$ is the degenerate distribution at $(x_0, \theta_0)$.
 The proof then follows essentially by obtaining the `drift' and `diffusion' coefficients of the discretized process, as in Equations (\ref{drift equations}) and then finding its limit. Formally, first obtain the quantities :
\beo
\mathbf{a}_n(\mathbf{y},t) &:=& \left( a_{n,i,j}(\my,t) \right)_{i,j=1,2} := n \int_{\mathbb{R}}(\mathbf{z}-\mathbf{y})(\mathbf{z}- \mathbf{y})' \Pi_n(\mathbf{y},d\mathbf{z}),\\
\mathbf{b}_n(\mathbf{y},t) &:=& \left( b_{n,k}(\my,t) \right)_{k=1,2} := n \int_{\mathbb{R}}(\mathbf{z}-\mathbf{y}) \Pi_n(\mathbf{y},d\mathbf{z})
.\eeo
The above is obtained by replacing $h^{-1}$ by $n$ in Equation (\ref{drift equations}).\\
Then find the matrix $\mathbf{a}$ and the vector 
$\mathbf{b}$ such that $\lim_{n \rightarrow \infty}|| \mathbf{a}_n(\my,t)-\mathbf{a}(\my,t)|| =0 $ and 
$\lim_{n \rightarrow \infty} ||\mathbf{b}_n(\my,t)-\mathbf{b}(\my,t)|| = 0$. Obtain the square root of matrix 
$\mathbf{a}(\mathbf{y},t)$(say $\mathbf{\sigma}_(\my,t)$), which satisfies 
$\mathbf{a}(\my,t) = \mathbf{\sigma}(\my,t) \mathbf{\sigma}(\my,t)^T$. 
These coefficients define a diffusion process uniquely which is non-explosive (see Remark \ref{conditions}), and the limiting process is governed by the 
equation: 
\beo
d \mathbf{Y}_t = \mathbf{b}(\mY_t,t) dt + \mathbf{\sigma}(\mY_t,t) d \mathbf{W}_t ,
\eeo
where $\mathbf{W}_t$ is a two dimensional Wiener process. 
For the processes defined in (\ref{X}) and (\ref{theta}), the limiting quantities $\mathbf{a}_n({\my,t})$ and $\mathbf{b}_n(\my,t)$ are \\
(for $\my=(x, \ \ \theta)$ ):
\beo
\lim_{n \rightarrow \infty} b_{n,1}(\my,t) &=& \frac{\theta^2}{2} \frac{\psi'(x)}{\psi(x)}, \\
\lim_{n \rightarrow \infty} b_{n,2}(\my,t) &=& \theta(q- \frac{\theta}{\sqrt{2\pi}}\frac{|\psi'(x)|}{\psi(x)}), \\
\lim_{n \rightarrow \infty} a_{n,1,1}(\my,t) &=& \theta^2, \\
\lim_{n \rightarrow \infty} a_{n,2,2}(\my,t) &=& 0, \\
\lim_{n \rightarrow \infty} a_{n,2,1}(\my,t) &=& 0  = 
\lim_{n \rightarrow \infty} a_{n,1,2}(\my,t)
\eeo
See Section \ref{drift and diffusion coefficients} for the derivations. \\
Since the trace norm of a matrix is a continuous function of its components we can say that 
$$||\mathbf{a}_n(\my,t)-\mathbf{a}(\my,t)|| \rightarrow 0  \ \mbox{and} \ ||\mathbf{b}_n(\my,t)-\mathbf{b}(\my,t)|| \rightarrow 0$$ 
where 
\beo
\mathbf{a}(\my,t) &=& \left(
\begin{array}{ll}
\theta^2 & 0\\
0 & 0
\end{array}
\right)
\Rightarrow \mathbf{\sigma}(\my,t)= \left(
\begin{array}{ll}
\theta & 0 \\
0 & 0
\end{array}
\right) \\
\mbox{and}\ \ \mathbf{b}(\my,t) &=& \left( \frac{\theta^2}{2} \frac{\psi'(x)}{\psi(x)},    \ \ \   \theta ( q - \frac{\theta}{\sqrt{2\pi}}\frac{|\psi'(x)|}{\psi(x)})  \right)^{T}
.\eeo
This proves the Theorem. 
\hfill{$\blacksquare$}

\section{Drift and diffusion coefficients}
\label{drift and diffusion coefficients}
Writing $\my=(x, \ \ \theta)$ we have 

\subsection{$b_{n,1}$}


\beo
&& b_{n,1}(\my,t) \\
&=& n E(X_n(\frac{i+1}{n})-X_n(\frac{i}{n})|\  \mY_n(\frac{i}{n}) = \my), \ \ \ \forall i=0,1,\ldots, \forall n \ge 1\\
&=& E(\sqrt{n} \theta_n(\frac{i}{n}) \xi_n(\frac{i+1}{n}) \epsilon_n(\frac{i+1}{n}) | \ \mY_n(\frac{i}{n}) = \my) \\
&=& \sqrt{n} \theta \Big( E(\xi_n(\fini) \epsilon_n(\fini) I_{A_n} |\  X_n(\fin)=x, \ \theta_n(\fin) = \theta)  \\
&+& E(\xi_n(\fini) \epsilon_n(\fini) I_{A_n^c} |\  X_n(\fin)=x, \ \theta_n(\fin) = \theta) \Big)
.\eeo
where $A_n(=A_n(x,\theta))$ is the set where $\xi_n(\fini)$ is one with probability 1, i.e, 
\beo
A_n(x, \theta) &=& \{ y: \frac{\psi(x +\frac{1}{\sqrt{n}} \theta y)}{\psi(x)} > 1 \}.\\
\mbox{Thus,} \ \ \ \lim_{n \rightarrow \infty}A_n^c(x, \theta) &=& \Big \{ 
\begin{array}{ll}
(-\infty,0)& \mbox{if $\psi'(x)>0$}\\
(0,\infty)& \mbox{if $\psi'(x)<0$} .
\end{array}
\eeo
 Therefore,
 \beo
 b_{n,1}(\my,t) &=& \sqrt{n} \theta \Big( \int_{A_n} \epsilon \phi(\epsilon) d \epsilon +  \int_{A_n^c} \frac{\psi(x + \frac{1}{\sqrt{n}} \theta \epsilon)}{\psi(x)} \epsilon \phi(\epsilon) d \epsilon \Big) \\
 &=& \sqrt{n} \theta \Big( \int_{A_n} \epsilon \phi(\epsilon) d \epsilon + \int_{A_n^c} \epsilon \phi(\epsilon) d \epsilon \\
 &+& \frac{\theta}{\sqrt{n}} \frac{\psi'(x)}{\psi(x)} \int_{A_n^c} \epsilon^2 \phi(\epsilon) d\epsilon 
 + O(\frac{1}{n})\Big),  \ \ \mbox{by Taylor's expansion,} \\
 &=& \sqrt{n} \theta \Big( \int_{\mathbb{R}} \epsilon \phi(\epsilon) d \epsilon +  \frac{\theta}{\sqrt{n}} \frac{\psi'(x)}{\psi(x)} \int_{A_n^c} \epsilon^2 \phi(\epsilon) d\epsilon + O(\frac{1}{n}) \Big) \\
 &=& \theta^2 \frac{\psi'(x)}{\psi(x)} \int_{A_n^c} \epsilon^2 \phi(\epsilon) d \epsilon + O(\frac{1}{\sqrt{n}}) \\
 \Rightarrow \lim_{n \rightarrow \infty} b_{n,1}(\my,t) &=& \theta^2 \frac{\psi'(x)}{\psi(x)} \lim_{n\rightarrow \infty} \int_{A_n^c} \epsilon^2 \phi(\epsilon) d \epsilon \n \\
&=& \Big \{ 
 \begin{array}{ll}
  \theta^2 \frac{\psi'(x)}{\psi(x)} \int_{-\infty}^{0} \epsilon^2 \phi(\epsilon) d \epsilon & \mbox {if $\psi'(x)>0$} \\
  \theta^2 \frac{\psi'(x)}{\psi(x)} \int_{0}^{\infty}\epsilon^2 \phi(\epsilon)d \epsilon & \mbox  {if $\psi'(x)<0$}
  \end{array} \\
  &=& \frac{\theta^2}{2} \frac{\psi'(x)}{\psi(x)} 
   .\eeo
 

 \subsection{$b_{n,2}$}
 \beo
b_{n,2}(\my,t) &=& n E(\theta_n(\frac{i+1}{n})- \theta_n(\fin)| \mY_n(\fin)=\my), \ \ \ \forall i=0,1,\ldots \\
&=& n E\Big( \theta_n(\fin) \{ e^{\frac{1}{\sqrt{n}}(\xi_n(\fini) - q_n(\fin))} -1 \} | \mY_n(\fin)= \my \Big) \\
&=& n \theta \Big( \frac{1}{\sqrt{n}} E(\xi_n(\fini) - q_n(\fin)|\mY_n(\fin)=\my) \\
&+& 
 E(\frac{1}{2n}(\xi_n(\fini) - q_n(\fin))^2|\mY_n(\fin)=\my) + O(\frac{1}{n^{3/2}}) \Big)\\
&=& \theta \sqrt{n} E(\xi_n(\fini) - q_n(\fin)|\mY_n(\fin)=\my) \\
&+& 
\frac{\theta}{2} E((\xi_n(\fini) - q_n(\fin))^2|\mY_n(\fin)=\my) 
+ O(\frac{1}{\sqrt{n}})
.\eeo
Now,
\be
& &\theta \sqrt{n} E(\xi_n(\fini) - q_n(\fin)|\mY_n(\fin)=\my) \n \\
&=& \theta \sqrt{n} \Big(E(\xi_n(\fini)| \mY_n(\fin)=\my) - q_n(\fin) \Big) \n \\
&=& \theta \sqrt{n} \Big( \int_{A_n} \phi(\epsilon) d \epsilon +
 \int_{A_n^c} \frac{\psi(x + \frac{1}{\sqrt{n}} \theta \epsilon)}{\psi(x)} \phi(\epsilon) d \epsilon - q_n(\fin)\Big) \n \\
&=& \theta \sqrt{n} \Big( \int_{A_n} \phi(\epsilon) d \epsilon \n \\
&+& \int_{A_n^c} \{ 1 +\frac{\theta}{\sqrt{n}} \frac{\psi'(x)}{\psi(x)} \epsilon + O(\frac{1}{n}) \} \phi(\epsilon) d \epsilon 
-  q_n(\fin) \Big) \n \\
&=&  \theta \sqrt{n} (1 - q_n(\fin)) \n \\
&+& \theta^2 \frac{\psi'(x)}{\psi(x)} \int_{A_n^c} \epsilon \phi(\epsilon) d \epsilon  + O(\frac{1}{\sqrt{n}}). \label{b_n,2,1}
\ee
And, 
\be
& & E\Big( (\xi_n(\fini) - q_n(\fin))^2| \mY_n(\fin)=\my \Big) \n \\
&=& E\Big( \xi_n(\fini)^2  | \mY_n(\fin)= \my \Big) \n \\
&-& 2q_n(\fin) E\Big(  \xi_n(\fini) | \mY_n(\fin)=\my \Big) +
q_n(\fin)^2  \n \\
&=& \int_{A_n} \phi(\epsilon) d \epsilon + \int_{A_n^c} \frac{\psi(x + \frac{1}{\sqrt{n}} \theta \epsilon)}{\psi(x)} \phi(\epsilon) d \epsilon  \n \\
&-& 2q_n(\fin) \Big(\int_{A_n} \phi(\epsilon) d \epsilon  + \int_{A_n^c} \frac{\psi(x + \frac{1}{\sqrt{n}} \theta \epsilon)}{\psi(x)} \phi(\epsilon) d\epsilon   \Big)  \n \\
&+& 
q_n(\fin)^2 \n \\
&=& (1-q_n(\fin))^2  + \frac{1}{\sqrt{n}} (1-2q_n(\fin)) \theta \frac{\psi'(x)}{\psi(x)} \int_{A_n^c} \epsilon \phi(\epsilon) d \epsilon \n \\
&+& O(\frac{1}{n}) 
\longrightarrow  0 , \label{b_n,2,2}
\ee
as $n \to \infty $ (since $ 1-q_n(\fin) \approx \frac{q}{\sqrt{n}}$), \tg{therefore}
\be 
\frac{1}{\sqrt{n}}(1-2q_n(\fin)) &\approx& \frac{1}{\sqrt{n}}(\frac{2q}{\sqrt{n}}-1). \n
\ee
Thus, from (\ref{b_n,2,1}) and (\ref{b_n,2,2}) we have, 
\beo
\lim_{n \rightarrow \infty} \mathbf{b}_{n,2}(\my,t) &=& \theta q + \theta^2 \frac{\psi'(x)}{\psi(x)} \lim_{n \rightarrow \infty} \int_{A_n^c} \epsilon \phi(\epsilon) d \epsilon \\
&=& \Big \{
\begin{array}{ll}
\theta \Big(q+\frac{\theta}{\sqrt{2\pi}} \frac{\psi'(x)}{\psi(x)}   \Big) & \mbox{if $\psi'(x)<0$}\\
\theta \Big(q-\frac{\theta}{\sqrt{2\pi}} \frac{\psi'(x)}{\psi(x)}   \Big) & \mbox{if $\psi'(x)>0$}
\end{array}
\\
&=& \theta \Big(q - \frac{\theta}{\sqrt{2\pi}} \frac{|\psi'(x)|}{\psi(x)} \Big) 
.\eeo


\subsection{$a_{n,1,1}$.}
\beo
a_{n,1,1}(\my,t) &=& n E \Big(  (X_n(\frac{i+1}{n})- X_n(\fin)^2) | \mY_{n}(\fin)= \my \Big) \ \, \forall i=0,1,\ldots\\
&=& \theta^2 E(\xi_n(\fini) \epsilon_n(\fini)^2|\mY_n(\fin)=\my)\\
&=& \theta^2 \Big( E(\xi_n(\fini) \epsilon_n(\fini)^2 I_{A_n}|\ \mY_n(\fin)=\my) \n \\
&+& 
E(\xi_n(\fini) \epsilon_n(\fini)^2I_{A_n^c}| \ \mY_n(\fin)=\my)  \Big)\\
&=& \theta^2 \Big(  \int_{A_n} \epsilon^2 \phi(\epsilon)d\epsilon + \int_{A_n^c}\epsilon^2 \frac{\psi(x+ \frac{1}{\sqrt{n}} \theta \epsilon)}{\psi(x)} \phi(\epsilon) d \epsilon \Big) \\
&=& \theta^2 \Big( \int_{A_n} \epsilon^2 \phi(\epsilon) d\epsilon + \int_{A_n^c}\epsilon^2 \phi(\epsilon) d\epsilon + O(\frac{1}{\sqrt{n}})  \Big)\\
&=& \theta^2 + O(\frac{1}{\sqrt{n}}). \\
\Rightarrow \lim_{n \rightarrow \infty} a_{n,1,1}(t) &=& \theta^2 .
\eeo
\subsection{$a_{n,2,2}$.}
\beo
a_{n,2,2}(\my,t) &=& n E\Big( (\theta_n(\frac{i+1}{n})-\theta_n(\frac{i}{n}))^2  | \mY_n(\fin)=\my \Big) \\
&=&n E\Big( \theta_n(\fin)^2 (e^{\frac{1}{\sqrt{n}} (\xi_n(\fini)-q_n(\fin))}-1)^2 |\mY_n(\fin) = \my \Big) \\ 
&=& n \theta^2 E \Big(  \Big\{ \frac{1}{\sqrt{n}} (\xi_n(\fini)-q_n(\fin)) + \frac{1}{2n} (\xi_n(\fini)\n \\
&-& q_n(\fin))^2 
+ O(\frac{1}{n^{3/2}}) \Big\}^2|\mY_n(\fin)=\my \Big)\\
&=& \theta^2 E\Big( (\xi_n(\fini)- q_n(\fin))^2 | \mY_n(\fin)=\my \Big) \n \\
&+& O(\frac{1}{\sqrt{n}})\\
\Rightarrow \lim_{n \rightarrow \infty} a_{n,2,2}(\my,t) &=& \theta^2 \lim_{n \rightarrow \infty} E\Big( \Big(\xi_n(\fini)-q_n(\fin)\Big)^2|\mY_n(\fin)=\my \Big) 
= 0,
\eeo
from (\ref{b_n,2,2}).

\subsection{$a_{n,1,2}$ and $a_{n,2,1}$.}
\beo
a_{n,1,2}(\my,t) 
&=& n E \Big(  \{X_n(\frac{i+1}{n})-X_n(\frac{i}{n})\} \{\theta_n(\frac{i+1}{n})- \theta_n(\fin)\} | \mY_n(\fin)=\my \Big)\\
&=& n E \Big(  \{ \fracn \theta_n(\fin) \xi_n(\fini) \epsilon_n(\fini) \} \{  \theta_n(\fin) (e^{\frac{1}{\sqrt{n}}(\xi_n(\fini)- q_n(\fin))}-1) \} \Big)\\
&=& \sqrt{n} \theta^2 E \Big(  \xi_n(\fini) \epsilon_n(\fini) \Big\{ \frac{1}{\sqrt{n}}(\xi_n(\fini)-q_n(\fin)) \n \\
&+& O(\frac{1}{n}) \Big\} | \mY_n(\fin)= \my \Big) \n \\
&=&  \theta^2 E\Big( \xi_n(\fini) \epsilon_n(\fini) (\xi_n(\fini)- q_n(\fin)) | \mY_n(\fin)=\my \Big) \\
&+& O(\frac{1}{\sqrt{n}}) . 
\eeo
Since \ $\xi_n = 0$, or $1$, \ $\xi_n^2 = \xi_n$. \ Hence \ $\xi_n\epsilon_n (\xi_n - q_n) = \xi_n^2 \epsilon_n - \xi_n \epsilon_n q_n = \xi_n \epsilon_n (1 - q_n)$.  
Therefore, 
\beo
&& E\Big( \xi_n(\fini) \epsilon_n(\fini) (\xi_n(\fini)-q_n(\frac{i}{n})) | \mY_n(\fin)= \my \Big)  \\
&=& (1-q_n(\frac{i}{n})) E\Big( \xi_n(\fini) \epsilon_n(\fini)  | \mY_n(\fin) = \my \Big)\\
&=& (1 - q_n(\frac{i}{n})) O(1) \  \longrightarrow 0 , \ \ \ \mbox{as } \ \ n \to \infty .
\eeo
Thus, 
$\lim_{n \rightarrow \infty} a_{n,1,2} =\lim_{n \rightarrow \infty} a_{n,2,1} = 0. \hfill{\blacksquare}$

\begin{rmk}
\label{mult ext}
{
Note that the form of the SDE in Theorem \ref{main thm 2} is similar to the Langevin diffusion equation (for univariate densities). 
This shows that this adaptive MCMC properly Normalized behaves in the limit as the Langevin diffusion which has 
$\psi(\cdot)$ as the invariant distribution. This bears a little resemblance to the Metropolis adjusted Langevin 
algorithm (MALA) procedure, where the proposal emulates the discretization of the Langevin algorithm. For more 
information regarding MALA and its convergence properties, see Marshall and Roberts \citep{Marshal}, Roberts and Rosenthal \citep{Roberts3}. }
\end{rmk}

\begin{rmk}
\label{conditions}
For a general target distribution $\psi(\cdot)$ we assume that the solutions satisfy the non-explosive condition given by, see Skorohod (\citep{Skorokhod})
\be
|b(\my,t)| + |\sigma(\my,t)| &\le& C (1 + |\my|),
\ee
for some constant $C>0$. We also assume that it also satisfies the local Lipschitz condition for uniqueness given by 
\be
|b(\my_1, t) - b(\my_2,t)| + ||\sigma(\my_1, t) - \sigma(\my_2, t)|| &\le& D_k(|\my_1 - \my_2|,
\ee
where $\my_1, \my_2 $ lies in some compact interval $S_k \subset \mathbb{R} \times \mathbb{R}^{+}$ and some constant $D_k>0$. Here 
\beo
||\sigma(\my, t)|| = \sqrt{\sum \limits \limits_{i,j=1}^{2} \sigma_{i,j}^2}. 
\eeo
For constant $\theta_t$ the non-explosive condition boils down to 
\be
\frac{|\psi'(x)|}{\psi(x)}&\le& C (1 + |x|), 
\ee
for some $C \ge 0.$
\end{rmk}

\begin{rmk}
\label{gronwall}
If the target density $\psi(\cdot)$ satisfy the linear growth condition
\beo
\frac{|\psi'(x)|}{\psi(x)}&\le& a |x| + b,
\eeo
 for some $a, b >0$, then from the SDE (\ref{coupled system}) we have that 
\beo
d \theta_s &=& \theta_s \Bl q - \fpi \frac{|\psi'(X_s)|}{\psi(X_s)}\Br ds \le q \theta_s ds \\
\R \theta_s &\le& \theta_0 e^{qs} \ \ \mbox{and,}\\
d X_s &=& \frac{\theta_s^2}{2} \frac{\psi'(X_s)}{\psi(X_s)} + \theta_s d W_s.
\eeo
Taking integrals from $0$ to $t$ we have 
\beo
X_t &=& X_0 + \intt \frac{\theta_s^2}{2}  \frac{\psi'(X_s)}{\psi(X_s)} ds  + \intt \theta_s  dW_s \\
\R |X_t| &\le& |X_0| + \intt \frac{\theta_s^2}{2}  \frac{|\psi'(X_s)|}{\psi(X_s)} ds +  |\intt \theta_s dW_s| \\
&\le& |X_0| + \intt \frac{\theta_s^2}{2}  \Bl a |X_s| + b \Br ds +  |\intt \theta_s dW_s| \\
&\le& |X_0| + \frac{a \theta_0^2}{2} \intt e^{2qs} |X_s| ds + \frac{b \theta_0^2}{4q}(e^{2qt} - 1)  + |\intt \theta_s dW_s|,
\eeo
using the bound for $\theta_t$. Taking expectations we have 
\beo
\R E(|X_t|) &\le& E(|X_0|) + \frac{a E(\theta_0^2)}{2} \intt e^{2qs} E(|X_s|) ds + \frac{b E(\theta_0^2)}{4q}(e^{2qt} - 1)  +  E(|\intt \theta_s dW_s|).
\eeo
By the Cauchy Schwarz inequality the last expectation is bounded by $$ \sqrt{E (\intt \theta_s d W_s)^2} = \sqrt{\intt E(\theta_s^2) ds} \le  \sqrt{E(\theta_0^2)} \sqrt{\intt e^{2qs} ds} = \frac{\sqrt{E(\theta_0^2)}}{\sqrt{2q}} \sqrt{e^{2qt} -1} \le  \frac{\sqrt{E(\theta_0^2)}}{\sqrt{2q}} e^{qt}.$$ Hence by a rearrangement of terms we have 
\beo
E|X_t|&\le& \underbrace{E(|X_0|)  + \frac{b E(\theta_0^2)}{4q}(e^{2qt} - 1)+ \frac{\sqrt{E(\theta_0^2)}} {\sqrt{2q}} e^{qt}}_{F_t}  + \intt \underbrace{\frac{a E(\theta_0^2)}{2}  e^{2qs}}_{A_s} E(|X_s|) ds \\
\eeo
Writing $G_t = E|X_t|, \ F_t = E(|X_0|)  + \frac{b E(\theta_0^2)}{4q}e^{2qt} + \frac{\sqrt{E(\theta_0^2)}}{\sqrt{2q}} e^{qt}, \ A_t = \frac{a E(\theta_0^2)e^{2qt}}{2}$ we have 
$$ G_t \le F_t + \int_{0}^{t} A_s G_s ds,$$
where $F_t$ is non negative and $A_t$ is increasing as a function of $t \in [0, \infty).$ Therefore from Gronwall's inequality, see, for example, \citep{Oksendal}, pp. 78, we have 
\beo
G_t &\le& F_t e^{\int_{0}^{t}A_s ds}, \ t \ge 0.
\eeo


Now $$
\intt A_s ds = \frac{a E(\theta_0^2)}{4q}(e^{2qt} -1) \le \frac{aE(\theta_0^2)}{4q} e^{2qt},$$ and so
\beo
E|X_t| &\le& \Bl E(|X_0|)  + \frac{b E(\theta_0^2)}{4q}e^{2qt} + \frac{\sqrt{E(\theta_0^2)}}{\sqrt{2q}} e^{qt} \Br e^{\frac{aE(\theta_0^2)}{4q} e^{2qt}}.
\eeo
This proves that the solution to the SDE of $(X_t,\theta_t$) given by Equation (\ref{coupled system}) is non-explosive. 
\end{rmk}

\section{Multi-dimensional target distribution}
\label{multsec}

In this section we consider the situation when the target distribution is a multivariate distribution $\psi(\mx), \ \mx \in \mR^p.$ Suppose the proposal distribution is multivariate Normal $N_p(\mathbf{0}, \Sigma)$. Then the adaptive algorithm will be given as:

\vskip 0.10in
\textbf{Algorithm 2} 
\benu
\item Select arbitrary $(\mX_0, \theta_0, \xi_0 ) \in \mR^{p} \times (0,\infty) \times \{0,1\}$. Set $n = 1$;
\item Propose a new move, say $\mY \sim N_p(\mX_{n-1}, \Sigma_{n-1})$ where $\Sigma_{n-1} = \theta_{n-1} \mathbf{I}_p,\ \mathbf{I}_p$ being the identity matrix of dimension $p$;
\item Accept the new point with probability $\alpha(\mX_{n-1}, \mY) = \min \{1, \frac{\psi(\mY)}{\psi(\mX_{n-1})}\}, \ \xi_n = 1$ if the sample is accepted else $\xi_n = 0$;
\item $ \theta_n = \theta_{n-1} e^{\frac{1}{\sqrt{n}}(\xi_n - q)}, \ q >0, \ \Leftrightarrow \log(\theta_n) = \log(\theta_{n-1}) + \frac{1}{\sqrt{n}} (\xi_n - q)$; 
\item $n \leftarrow n + 1 $ and go to step 2.
\eenu
This algorithm is equivalent to the following:\\

\textbf{Algorithm 2$'$:} \ \ \
\begin{enumerate}
\item Select $\{\mX_0,\theta_0,\xi_0 \}\in \mathbb{R}^p \times(0,\infty)\times\{0,1 \}, $
 where $\mathbb{R}^p$ is the state space. Set $n=1$;
\item Generate  $\bm{\epsilon}_{n-1} \sim N_p(\mathbf{0}, \Sigma_{n-1})$ where $\Sigma_{n-1} = \theta_{n-1} \mathbf{I}_p$. Given  $\mX_{n-1}, \theta_{n-1}, \bm{\epsilon_{n-1}} $ generate $$ \xi_n \sim \mbox{Bernoulli} \Bl \min \left \{  1, \frac{\psi(\mX_{n-1} + \theta_{n-1} \mathbf{\bm{\epsilon}_{n-1}}) }{\psi(\mX_{n-1})} \right \} \Br$$
and then
$$\mX_{n} = \mX_{n-1} + \theta_{n-1} \xi_{n} \bm{\epsilon}_{n-1}; $$ 
\item $\theta_{n}=\theta_{n-1} e^{\frac{1}{\sqrt{n}}(\xi_n - q)}, \ \ q>0,  \ \ \Leftrightarrow \log(\theta_n) =\log(\theta_{n-1})+ \frac{1}{\sqrt{n}}(\xi_n - q), \ \  q>0;$
\item
$n \leftarrow n+1$,
and go to step 2.
 \end{enumerate}
 
\begin{rmk}
 In Algorithm 2 all the co-ordinates $X_{in}, i=1,\ldots,p$, for a fixed $n \ge 1$, are scaled by the same factor $\theta_{n-1}$. This can be generalised where different co-ordinates are updated differently depending whether it is more mixing or not. We do not follow that approach here. 
 \end{rmk}
For the multivariate AMCMC, with the multivariate Normal proposal distribution, we now state the diffusion approximation  which is somewhat similar to the univariate AMCMC case as in Section \ref{diffusion approximation}. We also give the proof since it uses a slightly different method when compared to that of the univariate case and  requires the spherical symmetry property of the multivariate Normal $(\mathbf{0}, \mathbf{I}_p)$ distribution.
\begin{thm}
\label{multivariate lemma}
Applying the diffusion approximation (see Section 3 ) to Algorithm 2 such that $|| \nabla \psi(\mx)|| = 0$ on at most finitely many points, the diffusion corresponding to $\mY_t = (\mX_t, \theta_t)$ will be the solution of the following SDE:
\be
d \mY_t = \mb(\mY_t) dt + \mathbf{\sigma}(\mY_t) d \mathbf{W}_t, \label{mult sde}
\ee
where $\mb(\mY_t) = \Bl \frac{\theta_t^2}{2} \nabla \log \psi(\mX_t) , \  \theta_t(q - \fpi \theta_t ||\nabla \log \psi(\mX_t) ||) \Br^T,$ and 
\beo
\sigma(\mathbf{Y_t})=
\left(
\begin{array}{cc}
\theta_t \mathbf{I}_p & \mathbf{0}_{p \times 1} \\
\mathbf{0}_{1 \times p} & 0
\end{array}
\right).
\eeo
Here $\nabla \log \psi(\mX_t)= \frac{1}{\psi(\mX_t)} \Bl \frac{\partial}{\partial x_{1t}}  \psi(\mX_t), \frac{\partial}{\partial x_{2t}} \psi(\mX_t),\ldots, \frac{\partial}{\partial x_{pt}} \psi(\mX_t)\Br^{T} = \frac{\nabla \psi(\mX_t)^T}{\psi(\mX_t)}$ is the vector of partial derivatives of $\log \psi(\mx)$, $\mX_t = \Bl X_{1t}, X_{2t}, \ldots, X_{pt}\Br^{T}$ is the state vector, $\theta_t$ is the tuning parameter and $\bm{W}_t = \Bl W_{1t}, W_{2t}, \ldots, W_{(p+1)t}\Br^T$ is the $(p+1)$-dimensional Wiener process.\\
 \textbf{Proof:} Following the arguments and notations as in Sectiom \ref{diffusion approximation} we have to compute the `diffusion' and the `drift' coefficients which in this case are defined as: 
 \beo
 \mb_{n,1}(\my,t) &=& n E \Bl \mX_n(\frac{i+1}{n}) - \mX_n(\frac{i}{n})| \mY_n(\frac{i}{n}) = \my  \Br, \n \\
 \mb_{n,2}(\my,t) &=& n E \Bl \theta_n(\frac{i+1}{n}) - \theta_n(\frac{i}{n})| \mY_n(\frac{i}{n}) = \my  \Br, \n \\
 \mathbf{A}_{n,1,1}(\my, t)&=& n E \Bl (\mX_n(\frac{i+1}{n}) - \mX_n(\frac{i}{n})) (\mX_n(\frac{i+1}{n}) - \mX_n(\frac{i}{n}))^T| \mY_n(\frac{i}{n}) = \my \Br, \n \\
 \mathbf{A}_{n,2,2}(\my, t)&=& n E \Bl (\theta_n(\frac{i+1}{n}) - \theta_n(\frac{i}{n}))^2| \mY_n(\frac{i}{n}) = \my \Br, \ \mbox{and}  \n \\
 \mathbf{A}_{n,1,2}(\my, t)&=& n E \Bl (\theta_n(\frac{i+1}{n}) - \theta_n(\frac{i}{n}))(\mX_n(\frac{i+1}{n}) - \mX_n(\frac{i}{n}))| \mY_n(\frac{i}{n}) = \my \Br.
  \eeo
\end{thm}
Now, 
\beo
\mb_{n,1}(\my,t)&=& \sqrt{n} \theta \Big( E(\xi_n(\fini) \mathbf{\epsilon}_n(\fini) I_{A_n} |\  \mX_n(\fin)=\mx, \ \theta_n(\fin) = \theta)  \\
&+& E(\xi_n(\fini) \mathbf{\epsilon}_n(\fini) I_{A_n^c} |\  \mX_n(\fin)=\mx, \ \theta_n(\fin) = \theta) \Big)
\eeo
where $\mA _n (=\mA_n(\mx,\theta))$ is the set where $\xi_n(\fini)$ is one with probability 1, i.e, 
\beo
\mA_n &=& \{ \my: \frac{\psi(\mx +\frac{1}{\sqrt{n}} \theta \my)}{\psi(\mx)} \ge 1 \} \n \\
&=& \{ \my: (\psi(\mx) + \frac{1}{\sqrt{n}} \theta \nabla \psi(\mx)^T \my + O(\frac{1}{n}))/\psi(\mx) \ge 1\} \\
&=& \{ \my: \frac{1}{\sqrt{n}} \theta {\nabla \psi(\mx)}^T \my + O(\frac{1}{n}) \ge 0 \}.
\eeo
This implies that,
$$ \lim_{n \to \infty} \mA_{n} = \{\my: \nabla \psi(\mx)^T y  \ge 0	 \} := \mA \ (=\mA(\mx, \theta)). $$
 Therefore,
 \be
\mb_{n,1}(\my,t) &=& \sqrt{n} \theta \Big( \int_{\mA_n} \mathbf{\epsilon} \phi(\mathbf{\epsilon}) d \mathbf{\epsilon} +  \int_{\mA_n^c} \frac{\psi(\mx + \frac{1}{\sqrt{n}} \theta \mepsilon)}{\psi(\mx)} \mepsilon \phi(\mepsilon) d \mepsilon \Big)\n \\
&=& \sqrt{n} \theta \Bl \int_{\mA_n} \mepsilon \phi(\mepsilon) d \mepsilon + \int_{\mA_n^c} (1 + \frac{\theta}{\sqrt{n}\psi(\mx)} \nabla \psi(\mx)^T \mepsilon) \mepsilon \phi(\mepsilon) + O(\frac{1}{n})\Br \n \\
&=& \sqrt{n} \theta \int_{\mathbb{R}^p} \mepsilon \phi(\mepsilon) d\mepsilon  + \theta^2  \frac{1}{\psi(\mx)} \int_{\mA_n^c} ( \nabla \psi(\mx)^T \mepsilon) \mepsilon \phi(\mepsilon) d \mepsilon  + O(\frac{1}{\sqrt{n}})\n \\
&=& \theta^2  \frac{1}{\psi(\mx)} \int_{\mA_n^c} ( \nabla \psi(\mx)^T \mepsilon) \mepsilon \phi(\mepsilon) d \mepsilon  + O(\frac{1}{\sqrt{n}})\n \\
\R \lim_{n \to \infty} \mb_{n,1}(\my,t)&=& \theta^2  \frac{1}{\psi(\mx)} \lim_{n \to \infty} \int_{\mA_n^c} ( \nabla \psi(\mx)^T \mepsilon) \mepsilon \phi(\mepsilon) d \mepsilon \n \\
&=& \theta^2 \frac{1}{\psi(\mx)} \int_{\mA^c} ( \nabla \psi(\mx)^T \mepsilon) \ \mepsilon \phi(\mepsilon) d \mepsilon.
 \label{mult_b_n_1}
\ee
Consider the transformation 
\be 
\mepsilon :\to \mathbf{P}^T \epsilon  = \mZ = (Z_1, Z_2, \ldots, Z_p)^T, \label{trans}
\ee 
where $\mathbf{P}$ is an orthogonal matrix whose first column is $\frac{\nabla \psi(\mx)}{||\nabla \psi(\mx)||}:=P_1$, whenever $|| \nabla \psi(\mx)|| \neq 0$. Therefore $\nabla \psi(\mx)^T \epsilon = (|| \nabla \psi(\mx)||P_1)^T \epsilon  = {|| \nabla \psi(\mx)||P_1^T \epsilon} = || \nabla \psi(\mx)|| Z_1 $. Correspondingly $\mepsilon = \mathbf{P} \mZ$. The Jacobian of the transformation (\ref{trans}) is 1 and $Z_i, i = 1,2,\ldots,p$ are i.i.d $N(0,1)$, since $\epsilon_i, \ i = 1,2,\ldots,p$ are also i.i.d standard Normal. The integral in the RHS of Equation (\ref{mult_b_n_1}) is therefore 
\be
\int_{\mA^c} ( \nabla \psi(\mx)^T \mepsilon ) \mepsilon \phi(\mepsilon) d \mepsilon 
&=& || \nabla \psi(\mx) || \int_{ \{Z_1 < 0 \} }  Z_1 \mathbf{P} \mZ \phi(\mz) d \mz \n \\
&=& || \nabla \psi(\mx) || \int_{ \{Z_1 < 0 \} } Z_1 \sum \limits \limits_{i=1}^{p} P_i Z_i \phi(\mz) d \mz \n \\
&=& || \nabla \psi(\mx) || \Bl P_1 \int_{\{Z_1 <0 \}} Z_1^2 \phi(\mz) d\mz +  \sum \limits \limits_{i = 2}^{p} P_i \int_{ \{ Z_1 < 0 \} } Z_1 Z_i \phi(\mz) d \mz \n \Br \\
&=& || \nabla \psi(\mx) || P_1 \int_{ \{ Z_1 < 0 \} } Z_1^2 \phi(\mz) d \mz \n \  \mbox{(since $Z_i$'s are independent and $E(Z_i) = 0$)}\\
&=& || \nabla \psi(\mx) || \frac{P_1}{2}  \n \\
&=& \frac{1}{2} \nabla \psi(\mx), \ \mbox{since $P_1 = \frac{\nabla \psi(\mx)}{|| \nabla \psi(\mx) ||}$.} \n 
\ee
 Therefore 
\be
\mb_{1}(\my,t) &=& \frac{\theta^2}{2} \frac{\nabla \psi(\mx)} {\psi(\mx)} = \frac{\theta^2}{2} \nabla \log \psi(\mx). \label {mult_b_1}
\ee

For $b_{n,2}(\my,t)$ we have 
  \beo
b_{n,2}(\my,t) &=& n E(\theta_n(\frac{i+1}{n})- \theta_n(\fin)| \mY_n(\fin)=\my), \ \ \ \forall i=0,1,\ldots \\
&=& n E\Big( \theta_n(\fin) \{ e^{\frac{1}{\sqrt{n}}(\xi_n(\fini) - q_n(\fin))} -1 \} | \mY_n(\fin)= \my \Big) \\
&=& n \theta \Big( \frac{1}{\sqrt{n}} E(\xi_n(\fini) - q_n(\fin) ) |\mY_n(\fin)=\my) \\
&+& 
 E(\frac{1}{2n}(\xi_n(\fini) - q_n(\fin))^2|\mY_n(\fin)=\my) + O(\frac{1}{n^{3/2}}) \Big)\\
&=& \theta \sqrt{n} E(\xi_n(\fini) - q_n(\fin))|\mY_n(\fin)=\my) \\
&+& 
\frac{\theta}{2} E((\xi_n(\fini) - q_n(\fin))^2 |\mY_n(\fin)=\my) 
+ O(\frac{1}{\sqrt{n}})
.\eeo
Now for the first term ,
\be
& &\theta \sqrt{n} E(\xi_n(\fini) - q_n(\fin)|\mY_n(\fin)=\my) \n \\
&=& \theta \sqrt{n} \Big(E(\xi_n(\fini)| \mY_n(\fin)=\my) - q_n(\fin) \Big) \n \\
&=& \theta \sqrt{n} \Big( \int_{A_n} \phi(\mepsilon) d \mepsilon +
 \int_{A_n^c} \frac{\psi(\mx + \frac{1}{\sqrt{n}} \theta \mepsilon)}{\psi(\mx)} \phi(\mepsilon) d \mepsilon - q_n(\fin)\Big) \n \\
&=& \theta \sqrt{n} \Big( \int_{\mA_n} \phi(\mepsilon) d \mepsilon 
+ \int_{\mA_n^c} \{ 1 +\frac{\theta}{\sqrt{n}} \frac{\nabla \psi(\mx)^T}{\psi(\mx)} \epsilon + O(\frac{1}{n}) \} \phi(\mepsilon) d \mepsilon 
-  q_n(\fin) \Big) \n \\
&=&  \theta \sqrt{n} (1 - q_n(\fin)) 
+  \frac{\theta^2 }{\psi(\mx)} \int_{\mA_n^c} \nabla \psi(\mx)^T \mepsilon \phi(\mepsilon) d \mepsilon  + O(\frac{1}{\sqrt{n}}).\n 
\ee
And for the second term , 
\be
& & E\Big( (\xi_n(\fini) - q_n(\fin))^2| \mY_n(\fin)=\my \Big) \n \\
&=& E\Big( \xi_n(\fini)^2  | \mY_n(\fin)= \my \Big) \n \\
&-& 2q_n(\fin) E\Big(  \xi_n(\fini) | \mY_n(\fin)=\my \Big) +
q_n(\fin)^2  \n \\
&=& \int_{\mA_n} \phi(\mepsilon) d \mepsilon + \int_{\mA_n^c} \frac{\psi(\mx + \frac{1}{\sqrt{n}} \theta \mepsilon)}{\psi(\mx)} \phi(\mepsilon) d \mepsilon  \n \\
&-& 2q_n(\fin) \Big(\int_{\mA_n} \phi(\mepsilon) d \mepsilon  + \int_{\mA_n^c} \frac{\psi(\mx + \frac{1}{\sqrt{n}} \theta \mepsilon)}{\psi(\mx)} \phi(\mepsilon) d\mepsilon   \Big)  + q_n(\fin)^2 \n \\
&=& (1-q_n(\fin))^2  + \frac{1}{\sqrt{n}} (1-2q_n(\fin)) \theta \frac{1}{\psi(\mx)} \int_{A_n^c} \nabla \psi(\mx)^T \epsilon \phi(\epsilon) d \epsilon 
+ O(\frac{1}{n}) \n\\
&\to &  0 , \n 
\ee
as $n \to \infty $, since as before we assume that $ 1-q_n(\fin) \approx \frac{q}{\sqrt{n}}$. \tg{Therefore}
\be 
\frac{1}{\sqrt{n}}(1-2q_n(\fin)) &\approx& \frac{1}{\sqrt{n}}(\frac{2q}{\sqrt{n}}-1). \n 
\ee
Thus, only the first term contributes and we have 
\be
\lim_{n \rightarrow \infty} \mathbf{b}_{n,2}(\my,t) &=& \theta q + \frac{\theta^2}{\psi(\mx)} \lim_{n \rightarrow \infty} \int_{\mA_n^c} (\nabla \psi(\mx)^{T} \epsilon ) \phi(\epsilon) d \epsilon \n \\
&=& \theta q + \frac{\theta^2}{\psi(\mx)} \int_{\mA^c} ( \nabla \psi(\mx)^T \mepsilon ) \phi(\mepsilon) d \mepsilon. \label{multb2}
\ee
Using the transformation used in Equation (\ref{trans}) above we have 
\beo
\int_{\mA^c} (\nabla \psi(\mx)^T) \epsilon \phi(\mepsilon) d \mepsilon &=& || \nabla \psi(\mx)|| \int_{\{Z_1 < 0 \}} Z_1 \phi(\mz) d \mz \n \\
&=& || \nabla \psi(\mx)|| E (Z_1 I(Z_1 < 0)) \n \\
&=& -\fpi || \nabla \psi(\mx) ||.
\eeo
Therefore from (\ref{multb2}) we have
$$ b_2(\my, t) = \theta \Bl q - \fpi \frac{|| \nabla \psi(\mx)||}{\psi(\mx)} \theta \Br  = \theta \Bl q - \fpi \theta {|| \nabla \log \psi(\mx)||} \Br$$  

\beo
\mathbf{A}_{n,1,1}(\my,t) &=& n E \Big( (\mX_n(\frac{i+1}{n})- \mX_n(\fin))(\mX_n(\frac{i+1}{n})- \mX_n(\fin))^T | \mY_{n}(\fin)= \my \Big) \ \, \forall i=0,1,\ldots\\
&=& \theta^2 E(\xi_n(\fini) \mepsilon_n(\fini)\mepsilon_n(\fini)^T|\mY_n(\fin)=\my)\\
&=& \theta^2 \Big( E(\xi_n(\fini) \mepsilon_n(\fini)\mepsilon_n(\fini)^T I_{A_n}|\ \mY_n(\fin)=\my) \n \\
&+&E(\xi_n(\fini) \epsilon_n(\fini) \epsilon_n(\fini)^T I_{A_n^c}| \ \mY_n(\fin)=\my)  \Big)\\
&=& \theta^2 \Big(  \int_{A_n} \mepsilon \mepsilon^T \phi(\mepsilon)d\mepsilon + \int_{A_n^c}\mepsilon \mepsilon^T \frac{\psi(\mx+ \frac{1}{\sqrt{n}} \theta \epsilon)}{\psi(\mx)} \phi(\mepsilon) d \mepsilon \Big) \\
&=& \theta^2 \Big( \int_{A_n} \mepsilon\mepsilon^T \phi(\mepsilon) d\mepsilon + \int_{A_n^c} \epsilon \epsilon^T \phi(\mepsilon) d\mepsilon + O(\frac{1}{\sqrt{n}})  \Big)\\
&=& \theta^2 \int_{\mathbb{R}^p} \mepsilon \mepsilon^T \phi(\mepsilon) d \mepsilon + O(\frac{1}{\sqrt{n}}) = \theta^2 \mathbf{I}_p + O(\frac{1}{\sqrt{n}}). \\
\Rightarrow \lim_{n \to  \infty} \mathbf{A}_{n,1,1}(\my, t) &=& \theta^2  \mathbf{I}_p.
\eeo
The computations for $\mathbf{A}_{2,2}(\my,t)$ is same as that of the univariate case and is not repeated here.
\beo
\mathbf{A}_{n,1,2}(\my,t) 
&=& n E \Big(  \{\mX_n(\frac{i+1}{n})-\mX_n(\frac{i}{n})\} \{\theta_n(\frac{i+1}{n})- \theta_n(\fin)\} | \mY_n(\fin)=\my \Big)\\
&=& n E \Big(  \{ \fracn \theta_n(\fin) \xi_n(\fini) \mepsilon_n(\fini) \} \{  \theta_n(\fin) (e^{\frac{1}{\sqrt{n}}(\xi_n(\fini)- q_n(\fin))}-1) \} |\mY_n(
\frac{1}{n}) = \my \Big)\\
&=& \sqrt{n} \theta^2 E \Big(  \xi_n(\fini) \mepsilon_n(\fini) \Big\{ \frac{1}{\sqrt{n}}(\xi_n(\fini)-q_n(\fin)) + O(\frac{1}{n}) \Big\} | \mY_n(\fin)= \my \Big) \n \\
&=&  \theta^2 E\Big( \xi_n(\fini) \mepsilon_n(\fini) (\xi_n(\fini)- q_n(\fin)) | \mY_n(\fin)=\my \Big) \\
&+& O(\frac{1}{\sqrt{n}}) . 
\eeo
Since \ $\xi_n = 0$, or $1$, \ $\xi_n^2 = \xi_n$. \ Hence \ $\xi_n \mathbf{\epsilon_n} (\xi_n - q_n) = \xi_n^2 \mepsilon_n - \xi_n \mepsilon_n q_n = \xi_n \mepsilon_n (1 - q_n)$.  
Therefore, 
\beo
&& E\Big( \xi_n(\fini) \mepsilon_n(\fini) (\xi_n(\fini)-q_n(\frac{i}{n})) | \mY_n(\fin)= \my \Big)  \\
&=& (1-q_n(\frac{i}{n})) E\Big( \xi_n(\fini) \mepsilon_n(\fini)  | \mY_n(\fin) = \my \Big)\\
&=& (1 - q_n(\frac{i}{n})) O(1) \  \longrightarrow 0 , \ \ \ \mbox{as } \ \ n \to \infty .
\eeo
Similar computation for $\mathbf{A}_{n,2,1}$ yields that 
$\lim_{n \rightarrow \infty} \mathbf{A}_{n,2,1} =\lim_{n \rightarrow \infty} \mathbf{A}_{n,1,2} = 0$.
The drift and the diffusion coefficient gives the limiting diffusion of the multivariate AMCMC and this proves the theorem. $\hfill{\blacksquare}$

\begin{rmk}
For the uniqueness and the non-explosion of the solutions of the above SDE  (\ref{mult sde}) we need the local Lipschitz and the linear growth conditions as in Remark 2 and 3. This, in particular, would mean that $\nabla \log \psi(\cdot)$ satisfies the linear growth condition
\be
||\nabla \log \psi(\mx)|| &\le& a ||\mx|| + b, \ \forall \mx \in \mathbb{R}^p, \label{mult_growth}
\ee
for some $a>0$ and $b \ge 0$. 
\end{rmk}



\section{Simulation}
\label{sim}
We give some plots of the Adaptive and the Non Adaptive samplers with Normal(0,1) as the target density for different choice of starting $\theta_0$ and $q = 0.50$ (see Figure 1 and 2). For the Non Adaptive chain the $\theta_n$ is kept constant at $\theta_0$. 
The number of samples used as burn-in was 1000 in both the cases. The density function of $N(0,1)$ is overlapped on the histogram of the sample generated by Standard MCMC and Adaptive MCMC for $\theta = 1$ and $\theta = 10$.
Comparing the two figures it seems that adaptive MCMC is better at sampling form the target density $N(0,1)$ starting with a large value of $\theta_0.$
The plot of the sample generated from the diffusion corresponding to the Adaptive and Standard MCMC given by Equation (\ref{coupled system}) and Equation (\ref{X MC}), obtained by the Euler method, is given in Figure 3 


%

\begin{figure}[h]
\includegraphics[width=15cm, height=20cm]{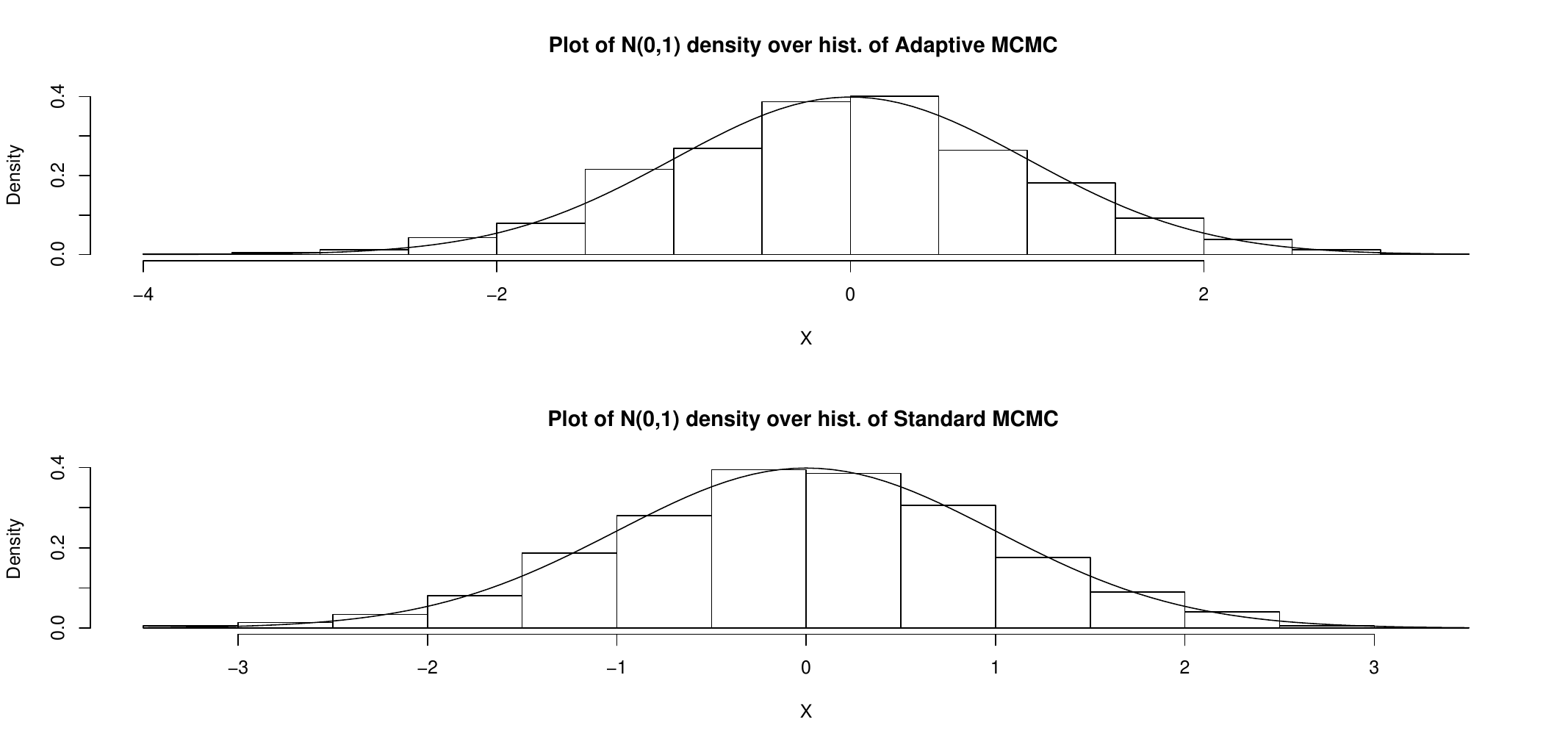}
\caption{Density of $N(0,1)$ overlapped on the histogram of the sample generated using Standard and Adaptive MCMC using $q = 0.50$ and $\theta_0=1$.}
\end{figure}

\begin{figure}[h]
\includegraphics[width=15cm, height=20cm]{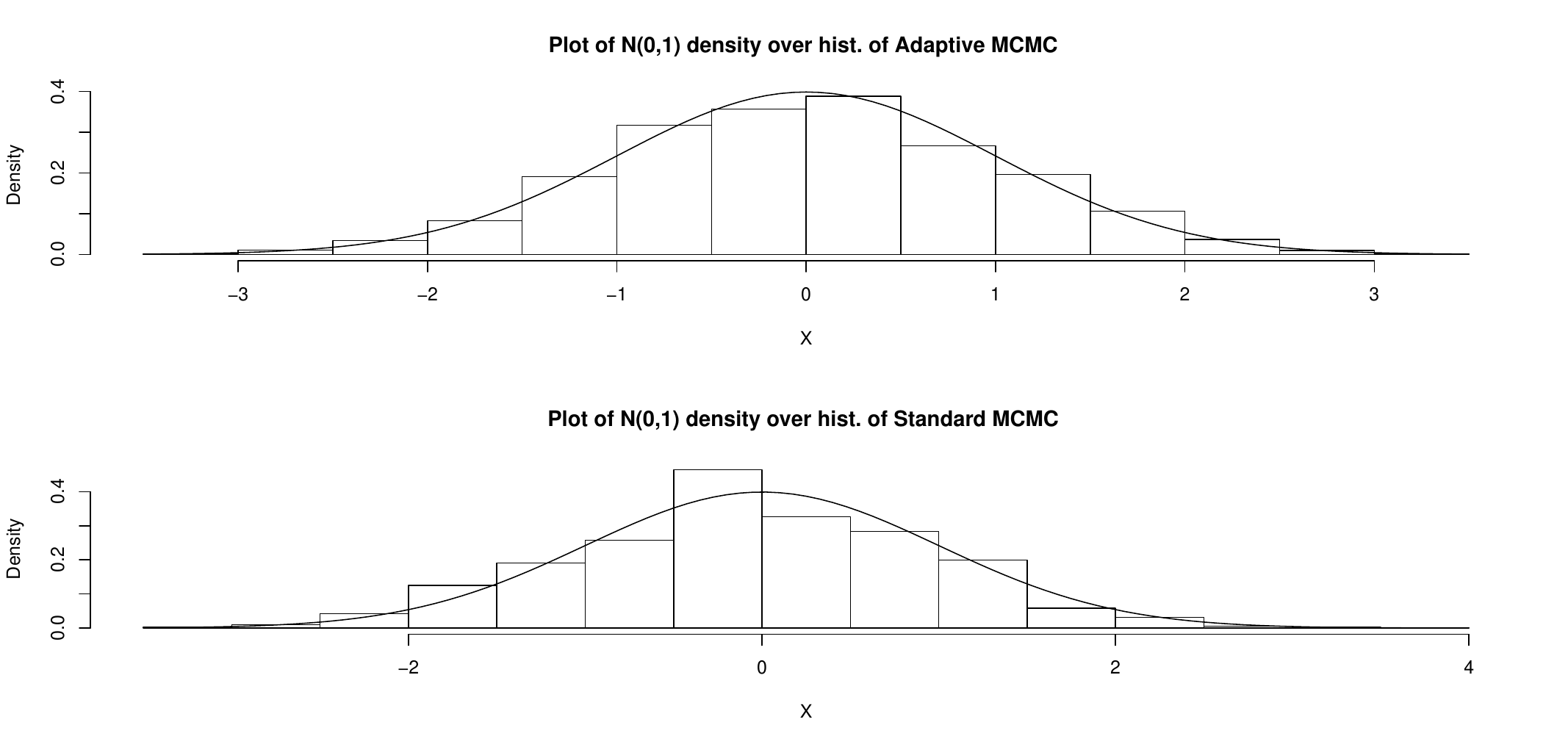}
\caption{Density of $N(0,1)$ overlapped on the histogram of the sample generated using Standard and Adaptive MCMC using $q = 0.50$ and $\theta_0=10$.}
\end{figure}

\begin{figure}[h]
\includegraphics[width=15cm, height=20cm]{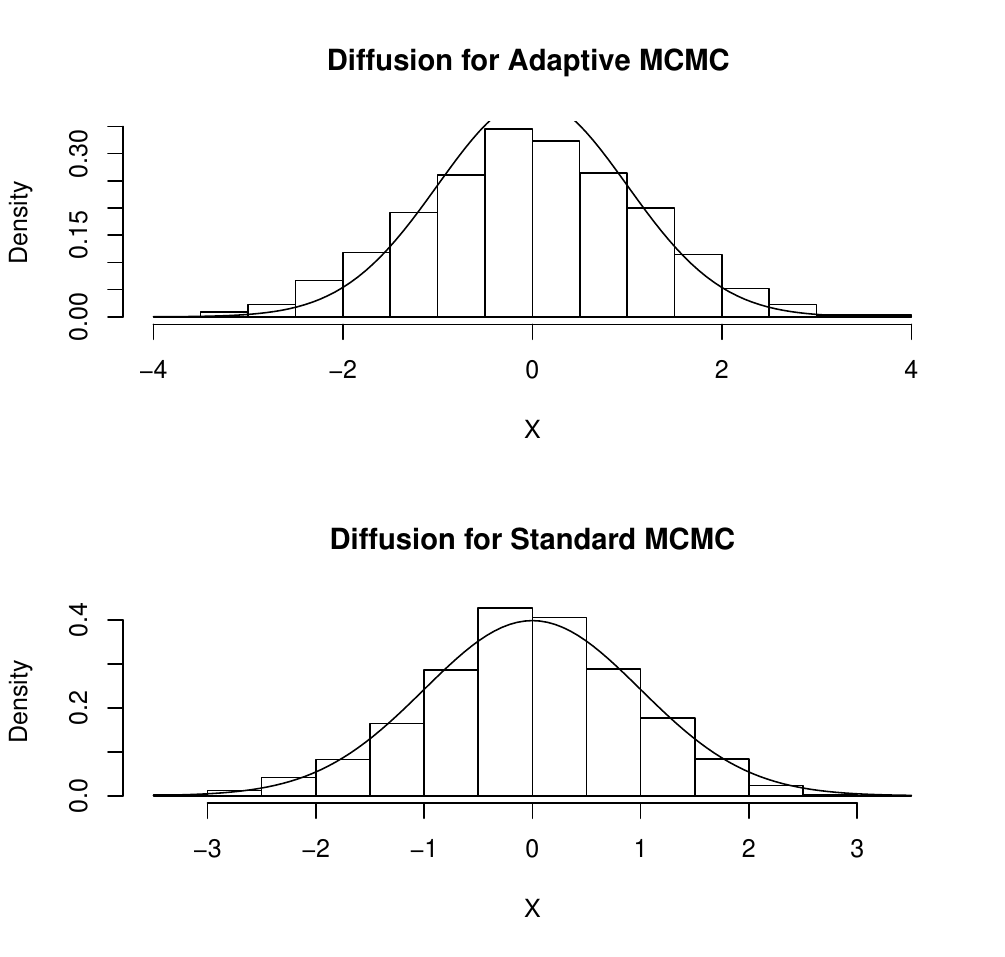}
\caption{Density of $N(0,1)$ overlapped on the histogram of the sample generated from the diffusion corresponding to Standard and Adaptive MCMC using $q = 1.0$ and $\theta_0=1$.}
\end{figure}


%
%


\section {Conclusion}
\label{summary}

Diffusion approximation is a well studied technique that has been applied to many fields (e.g., \citep{Ethier}, \citep{Nelson}). In AMCMC the tuning parameter changes as the iteration progresses and therefore the transition 
kernel also changes. As a result the invariant properties of the chain are not easily obtainable. In this paper we have applied the diffusion approximation 
procedure to the AMCMC for both univariate and multivariate target distributon using the standard univariate Normal and standard multivariate Normal distribution as the propsal.  In both the cases we obtain the limiting diffusion. Although the procedure can be extended to any univariate proposal with finite second moments and symmetric about 0, such extension to the multivariate proposal is not straight forward. 
{Diffusive limits for Metropolis Hastings 
 algorithm were earlier obtained in \citep{Roberts3,Stramer1, Stramer2}. 
Also, there are some recent work on diffusive limits of high-dimensional non-adaptive MCMC, see, for example, Mattingly \textit{et al}. \citep{Mattingly}).
Our technique expands the scope of comparison between AMCMC and Standard MCMC, as embedding in continuous time allows 
various discrete approximations through which one can compare them in finer details.\\


}




\end{document}